\documentstyle[12pt,amssym]{article}
\newcommand{\dotimes}{\stackrel{}{\dot{\otimes}}}

\def\cX{\mbox{$\cal X$}}
\def\cR{\mbox{$\cal R$}}
\def\id{\mbox{\rm id}}
\def\ts{\textstyle}
\def\half{\frac{1}{2}}
\def\thalf{\frac{3}{2}}
\def\tT{\tilde{T}}
\def\hJ{\hat{J}}
\def\N{\mbox{$\Bbb N$}}
\def\tR{\tilde{R}}
\def\hF{\hat{F}}
\def\hZ{\hat{Z}}
\def\ha{\hat{a}}
\def\hb{\hat{b}}
\def\hc{\hat{c}}
\def\hd{\hat{d}}
\def\hT{\hat{T}}
\def\hD{\hat{D}}
\def\ta{\tilde{a}}
\def\tb{\tilde{b}}
\def\tc{\tilde{c}}
\def\td{\tilde{d}}
\def\tD{\tilde{D}}
\def\cT{\mbox{$\cal T$}}
\def\halpha{\hat{\alpha}}
\def\hbeta{\hat{\beta}}
\def\hgamma{\hat{\gamma}}
\def\hdelta{\hat{\delta}}
\def\hcJ{\mbox{\hskip 6pt $\hat{\mbox{\hskip -6pt $\cal J$}}$}}

\title{
ON JORDANIAN U$_{h,\alpha}$(gl(2)) ALGEBRA AND ITS $T$
MATRICES VIA A CONTRACTION METHOD}
\author{R. CHAKRABARTI$^1$ and C. QUESNE$^{2,}$\thanks{Directeur de recherches
FNRS; E-mail address: cquesne@ulb.ac.be}\\ {\small $^1$ \sl Department of
Theoretical Physics, University of Madras,}\\ {\small \sl Guindy Campus,
Madras-600025, India}\\ {\small $^2$ \sl Physique Nucl\'eaire Th\'eorique et
Physique Math\'ematique,}\\ {\small \sl Universit\'e Libre de Bruxelles,
Campus de
la Plaine CP229,} \\ {\small \sl  Boulevard~du Triomphe, B-1050 Brussels,
Belgium}}
\date{ }
\begin{document}
\baselineskip=22pt plus 1pt minus 1pt
\maketitle

\begin{abstract}
The $R_h^{j_1;j_2}$ matrices of the Jordanian U$_h$(sl(2))
algebra at arbitrary dimensions may be obtained from the corresponding
$R_q^{j_1;j_2}$ matrices of the standard $q$-deformed U$_q$(sl(2)) algebra
through a contraction technique. By extending this method, the coloured
two-parametric ($h, \alpha$) Jordanian
$R_{h,\alpha}^{j_1,z_1;j_2,z_2}$ matrices of the U$_{h,\alpha}$(gl(2))
algebra may
be derived from the corresponding coloured
$R_{q,\lambda}^{j_1,z_1;j_2,z_2}$ matrices of the standard ($q,
\lambda$)-deformed U$_{q,\lambda}$(gl(2)) algebra. Moreover, by using the
contraction process as a tool, the coloured $T_{h,\alpha}^{j,z}$ matrices for
arbitrary ($j, z$) representations of the Jordanian Fun$_{h,\alpha}$(GL(2))
algebra
may be extracted from the corresponding $T_{q,\lambda}^{j,z}$ matrices of the
standard Fun$_{q,\lambda}$(GL(2)) algebra.
\end{abstract}

\newpage
%
%
\section{Introduction}
\label{sec:intro}
The Lie group GL(2) only admits (up to isomorphism) two distinct quantum group
deformations with central quantum determinant~\cite{kuper}: the standard
$q$-deformation Fun$_q$(GL(2))~\cite{drinfeld87}, and the Jordanian deformation
Fun$_h$(GL(2))~\cite{demidov}. On the quantum algebra level, the Jordanian
deformation U$_h$(sl(2)) of the classical enveloping algebra U(sl(2)) was first
considered in Ref.~\cite{ohn}, and its universal $\cR_h$ matrix was given in
Ref.~\cite{ballesteros}. The fundamental representation of U$_h$(sl(2)), which
remains undeformed, was obtained in Ref.~\cite{ohn}, while the other
finite-dimensional highest-weight representations were first studied in
Ref.~\cite{dobrev}.\par
%
%
Both Fun$_q$(GL(2)) and Fun$_h$(GL(2)) quantum groups may be viewed as special
cases of two-parametric deformations,
Fun$_{q,\lambda}$(GL(2))~\cite{schirrmacher} and
Fun$_{h,\alpha}$(GL(2))~\cite{agha93}, respectively. The quantum algebra
U$_{h,\alpha}$(gl(2)), dual to the two-parametric Jordanian quantum group
Fun$_{h,\alpha}$(GL(2)), was found in Ref.~\cite{aneva}. It was shown there that
with an appropriate choice of parameters, the algebraic properties of the
two-parametric Jordanian U$_{h,\alpha}$(gl(2)) algebra coincide with those
of the
single-parametric U$_h$(gl(2)) one. The second parameter~$\alpha$, however,
plays
a manifest role in the coalgebraic properties. An independent construction
of the
Hopf algebra U$_{h,\alpha}$(gl(2)) was carried out in Ref.~\cite{parashar},
together
with a determination of the corresponding universal $\cR_{h,\alpha}$ matrix.\par
%
%
Some related coloured Hopf algebraic structures were also introduced. In
Ref.~\cite{cq}, the single-parametric Jordanian U$_h$(gl(2)) algebra was
extended
to a coloured Hopf algebra U$^c$(gl(2)), endowed with a coloured universal
$\cR^c$ matrix. More recently, a coloured two-parametric Jordanian
$R_{h,\alpha}^{z_1;z_2}$ matrix, corresponding to the fundamental
representation,
was constructed, and used to build a coloured quantum group
Fun$_{h,\alpha}^{z_1;z_2}$(GL(2))~\cite{parashar}, thereby generalizing to the
Jordanian context a previous work in the standard $q$-deformed
context~\cite{basu}.\par
%
%
Two useful tools have been devised for studying Jordanian deformations. One of
them is a contraction procedure that allows one to construct the latter from
standard deformations~\cite{agha95}: a similarity transformation of the defining
$R_q$ and~$T_q$ (resp.~$R_{q,\lambda}$ and~$T_{q,\lambda}$) matrices of
Fun$_q$(SL(2)) (resp.~Fun$_{q,\lambda}$(GL(2))) is performed using a matrix
singular itself in the $q\to 1$ limit, but in such a way that the transformed
matrices are non-singular, and yield the defining $R_h$ and~$T_h$
(resp.~$R_{h,\alpha}$ and~$T_{h,\alpha}$) matrices of Fun$_h$(SL(2))
(resp.~Fun$_{h,\alpha}$(GL(2))). Such a contraction technique can be
generalized to
higher-dimensional quantum groups~\cite{ali}.\par
%
%
The other tool is a nonlinear invertible map between the generators of
U$_h$(sl(2)) and U(sl(2))~\cite{abdesselam96}. This map yields an explicit and
simple method for constructing the finite-dimensional irreducible
representations
of U$_h$(sl(2))~\cite{abdesselam96}, and for studying various related
problems~\cite{aizawa,joris}.\par
%
%
Recently, an operational generalization of the contraction method combining both
tools was proposed~\cite{abdesselam98}. It yields the Jordanian deformation of
U(sl($N$)) from the standard one along with a nonlinear map from the $h$-Borel
subalgebra on the corresponding classical Borel subalgebra, which in the case of
U$_h$(sl(2)) can be extended to the whole algebra. By using such a generalization,
the $R_h^{\half;j}$ matrix of the $\left(\half \otimes j\right)$ representation
of~U$_h$(sl(2)) was obtained from the corresponding $R_q^{\half;j}$ matrix
of the
standard $q$-deformed U$_q$(sl(2)) algebra~\cite{abdesselam98}. Furthermore, it
was shown that the Drinfeld twist operator~\cite{drinfeld90}, relating the usual
classical U(sl(2)) coalgebraic properties with the nontrivial coalgebraic
structure
of the Jordanian U$_h$(sl(2)) algebra, can be found as a series expansion in the
deformation parameter $h$~\cite{abdesselam98}.\par
%
%
The purpose of the present paper is threefold. Firstly, in
Sec.~\ref{sec:R_h}, we
wish to complete the work of Ref.~\cite{abdesselam98} by providing higher order
terms in the expansion of the Drinfeld twist operator, and by demonstrating that
the $h$-Jordanian $R_h^{j_1;j_2}$ matrices of arbitrary $(j_1 \otimes j_2)$
representations of U$_h$(sl(2)) may be obtained by contracting the corresponding
standard $R_q^{j_1;j_2}$ matrices.\par
%
%
Secondly, in Sec.~\ref{sec:R_{h,alpha}}, we wish to extend
Ref.~\cite{abdesselam98} in two ways: by going from the single-parametric
case to
the two-parametric one, and from $(j_1 \otimes j_2)$ standard representations to
$(j_1, z_1 \otimes j_2, z_2)$ coloured representations, where the colour
parameter~$z$ assumes arbitrary distinct values in different sectors of the
tensor
product space. We will prove that the contraction scheme provides a
mechanism for
extracting the $R_{h,\alpha}^{j_1,z_1;j_2,z_2}$ matrix in an arbitrary
$(j_1, z_1 \otimes j_2, z_2)$ coloured representation of the
U$_{h,\alpha}$(gl(2))
algebra from the corresponding $R_{q,\lambda}^{j_1,z_1;j_2,z_2}$ matrix of the
two-parametric standard U$_{q,\lambda}$(gl(2)) algebra. As a result, we will get
solutions of the coloured Yang-Baxter equation
\begin{equation}
  R_{12}^{j_1,z_1;j_2,z_2} R_{13}^{j_1,z_1;j_3,z_3} R_{23}^{j_2,z_2;j_3,z_3}
  = R_{23}^{j_2,z_2;j_3,z_3} R_{13}^{j_1,z_1;j_3,z_3} R_{12}^{j_1,z_1;j_2,z_2},
  \label{eq:YB}
\end{equation}
where we have suppressed the subscripts $(h, \alpha)$.\par
%
%
Thirdly, in Sec.~\ref{sec:T_{h,alpha}}, we wish to apply the contraction
process of
Ref.~\cite{abdesselam98} to the construction of an entirely new type of
Jordanian
deformed objects, namely the finite-dimensional coloured representations of $(h,
\alpha)$-deformed group-like elements $T_{h,\alpha}^{j,z}$ of
Fun$_{h,\alpha}$(GL(2)). In this respect, it should be stressed that if the
universal
$\cT_{q,\lambda}$ matrix, acting as a dual form between the $(q,
\lambda)$-deformed standard Hopf algebras Fun$_{q,\lambda}$(GL(2)) and
U$_{q,\lambda}$(gl(2)), is well known~\cite{fronsdal}, as well as its
finite-dimensional coloured representations $T_{q,\lambda}^{j,z}$~\cite{jaga},
nothing similar is available yet for Jordanian Hopf algebras (but for the
defining
$T_{h,\alpha}$ matrix). Here, we will prove that the $T_{q,\lambda}^{j,z}$
matrices
generate, through the contraction procedure, the corresponding Jordanian
$T_{h,\alpha}^{j,z}$ matrices. In this way, we will construct solutions of the
coloured inverse scattering equation
\begin{equation}
  R_{h,\alpha}^{j_1,z_1;j_2,z_2} \left(T_{h,\alpha}^{j_1,z_1} \otimes 1\right)
  \left(1 \otimes T_{h,\alpha}^{j_2,z_2}\right) =
  \left(1 \otimes T_{h,\alpha}^{j_2,z_2}\right)
  \left(T_{h,\alpha}^{j_1,z_1} \otimes 1\right) R_{h,\alpha}^{j_1,z_1;j_2,z_2}.
  \label{eq:RTT}
\end{equation}
%
%
\section{The $R_h$ Matrices of the U$_h$(sl(2)) Algebra as Contraction Limits}
\label{sec:R_h}
\setcounter{equation}{0}
We start by enlisting the Hopf structure of the single-parametric Jordanian
algebra U$_h$(sl(2))~\cite{ohn}. The algebra reads
\begin{equation}
  [H, X] = 2 \frac{\sinh hX}{h}, \quad [H, Y] = - Y (\cosh hX) - (\cosh hX) Y, \quad [X,
  Y] = H.  \label{eq:h-alg}
\end{equation}
The non-cocommutative coproduct, counit and the antipode maps assume the form
\begin{eqnarray}
  \Delta_h(X) & = & X \otimes 1 + 1 \otimes X, \qquad \Delta_h(Y) = Y
\otimes T +
           T^{-1} \otimes Y, \nonumber \\
  \Delta_h(H) & = & H \otimes T + T^{-1} \otimes H, \nonumber \\
  \epsilon_h(\cX) & = & 0 \qquad \forall \cX \in \{X, Y, H\}, \nonumber \\
  S_h(X) & = & - X, \qquad S_h(Y) = - TYT^{-1}, \qquad S_h(H) = - THT^{-1},
           \label{eq:h-coalg}
\end{eqnarray}
where $T = \exp(hX)$. The universal $\cR_h$ matrix of the triangular Hopf
algebra
U$_h$(sl(2)) is given~\cite{ballesteros} in a convenient form by
\begin{equation}
  \cR_h = \exp(- hX \otimes TH) \exp(hTH \otimes X).  \label{eq:univ-R_h}
\end{equation}
{}For a quasitriangular Hopf algebra~$U$, the universal $\cR$ matrix, an
invertible
element in $U \otimes U$, satisfies the relations
\begin{eqnarray}
  (\Delta \otimes \id) \cR & = & \cR_{13} \cR_{23}, \qquad (\id \otimes \cR) =
            \cR_{13} \cR_{12}, \nonumber \\
  \sigma \circ \Delta(\cX) & = & \cR \Delta(\cX) \cR^{-1} \qquad \forall
\cX \in U,
            \label{eq:R-prop}
\end{eqnarray}
where $\sigma$ is the permutation in $U \otimes U$. For a triangular Hopf
algebra, the additional property $\cR^{-1} = \cR_{21}$ holds.\par
%
%
An  invertible nonlinear map of the generating elements of the U$_h$(sl(2))
algebra on the elements of the classical U(sl(2)) algebra
plays~\cite{abdesselam98} a pivotal role in obtaining the $h$-Jordanian
$R_h^{\half;j}$ matrix in the $({\ts \half} \otimes j)$ representation as a
suitable
contraction limit of the corresponding standard $R_q^{\half;j}$ matrix. The map
reads~\cite{abdesselam98}
\begin{equation}
  T = \tT, \qquad Y = J_- - {\ts \frac{1}{4}} h^2 J_+ \left(J_0^2 -
1\right), \qquad
  H = \left(1 + (hJ_+)^2\right)^{1/2} J_0,  \label{eq:nl-map}
\end{equation}
where $\tT = h J_+ + \left(1 + (hJ_+)^2\right)^{1/2}$. The elements $(J_{\pm},
J_0)$ are the generators of the classical U(sl(2)) algebra:
\begin{eqnarray}
  [J_0, J_{\pm}] & = & \pm 2 J_{\pm}, \quad [J_+, J_-] = J_0,
\label{eq:class-alg} \\
  \Delta_0(J_{\zeta}) & = & J_{\zeta} \otimes 1 + 1 \otimes J_{\zeta}, \quad
          \epsilon_0(J_{\zeta}) = 0, \quad S_0(J_{\zeta}) = - J_{\zeta} \quad
          \forall \zeta \in \{\pm, 0\}.  \label{eq:class-coalg}
\end{eqnarray}
\par
%
%
{}Following Drinfeld's arguments~\cite{drinfeld90}, it is possible to
construct a
twist operator $G \in \mbox{U(sl(2))}^{{}\otimes 2}[[h]]$ relating the
$h$-Jordanian coalgebraic structure~(\ref{eq:h-coalg}) with the corresponding
classical structure~(\ref{eq:class-coalg}). For an invertible map $m: (X,
Y, H) \to
(J_{\pm}, J_0)$, $m^{-1}: (J_{\pm}, J_0) \to (X, Y, H)$, the following
relations hold:
\begin{eqnarray}
  (m \otimes m) \circ \Delta_h \circ m^{-1}(\cX) & = & G \Delta_0(\cX) G^{-1},
           \nonumber \\
  m \circ S_h \circ m^{-1}(\cX) & = & g S_0(\cX) g^{-1},
\end{eqnarray}
where $\cX \in \mbox{U(sl(2))}[[h]]$. The transforming operator~$g$ ($\in
\mbox{U(sl(2))}[[h]]$), and its inverse may be expressed as
\begin{equation}
  g = \mu \circ (\id \otimes S_0) G, \qquad g^{-1} = \mu \circ (S_0 \otimes \id)
  G^{-1},
\end{equation}
where $\mu$ is the multiplication map.\par
%
%
{}For the map~(\ref{eq:nl-map}), we have the following construction:
\begin{eqnarray}
  G & = & 1 \otimes 1 - {\ts \half} h r + {\ts \frac{1}{8}} h^2 \left[r^2 +
2 (J_+
          \otimes J_+) \Delta_0(J_0)\right] \nonumber \\
  & & \mbox{} - {\ts \frac{1}{48}} h^3 \left[r^3 + 6 (J_+ \otimes J_+)
\Delta_0(J_0)
          r - 4 \bigl(\Delta_0(J_+)\bigr)^2 r \right] \nonumber \\
  & & \mbox{} + {\ts \frac{1}{384}} h^4 \Bigl[r^4 - 16
\bigl(\Delta_0(J_+)\bigr)^2 r^2
          + 12 (J_+ \otimes J_+) \Delta_0(J_0) r^2 \nonumber \\
  & & \mbox{} + 12 \bigl((J_+ \otimes J_+) \Delta_0(J_0)\bigr)^2 + 6 \left(J_+^2
           \otimes 1 - 1 \otimes J_+^2\right)^2 \Delta_0(J_0) \nonumber \\
  & & \mbox{} + 12 \bigl(\Delta_0(J_+)\bigr)^2 \left(J_+^2 \otimes 1 + 1 \otimes
          J_+^2\right) \Delta_0(J_0) \nonumber \\
  & & \mbox{} - 8 \Delta_0(J_+) \left(J_+^3 \otimes 1 + 1 \otimes J_+^3\right)
          \Delta_0(J_0) - 10 \bigl(\Delta_0(J_+)\bigr)^4 \Delta_0(J_0)\Bigr]
          \nonumber \\
  & & \mbox{} + O\left(h^5\right), \nonumber \\
  g & = & 1 + hJ_+\bigl(1+h^2J_+^2 \bigr)^{1/2} + h^2J_+^2, \nonumber \\
  g^{-1} & = & 1 - hJ_+ \bigl( 1+h^2J_+^2 \bigr)^{-1/2},
          \label{eq:G-g}
\end{eqnarray}
where the Jordanian classical $r$ matrix reads $r = J_0 \otimes J_+ - J_+
\otimes J_0$. The transforming operator~$g$ is obtained in a closed form
in~(\ref{eq:G-g}).  The first non-trivial term $O\left(h^2\right)$ in the series
expansion of~$G$ in~(\ref{eq:G-g}) was previously evaluated~\cite{abdesselam98}.
The series expansion of the twist operator $G$ may be developed up to an
arbitrary order in~$h$. It may be verified that the expansion~(\ref{eq:G-g}), in
powers of~$h$, of the twist operator~$G$ corresponding to the
map~(\ref{eq:nl-map}) satisfies the cocycle condition
\begin{equation}
  (1 \otimes G) (\id \otimes \Delta_0) G = (G \otimes 1) (\Delta_0 \otimes
\id) G
\end{equation}
up to the desired order. The universal \cR$_h$ matrix~(\ref{eq:univ-R_h}),
used in
conjunction with the map~(\ref{eq:nl-map}), may be recast in the form
\begin{equation}
  \cR_h = (\sigma \circ G) G^{-1},
\end{equation}
valid up to an arbitrary order in the expansion~(\ref{eq:G-g}).\par
%
%
In Ref.~\cite{abdesselam98}, the $R_q^{\half;j}$ matrix of the standard
$q$-deformed U$_q$(sl(2)) algebra has been found to yield, after a suitable
transformation and a subsequent $q \to 1$ limiting process, the corresponding
$h$-Jordanian $R_h^{\half;j}$ matrix. A well-defined construction along
this route
for the $R_h^{\half;j}$ matrix exists as the singularities systematically
cancel.
Similar construction of $R_h^{j_1;j_2}$ matrices may be continued for arbitrary
$(j_1 \otimes j_2)$ representations.\par
%
%
To this end, we proceed with the universal \cR$_q$ matrix of the U$_q$(sl(2))
algebra, given by~\cite{majid}
\begin{equation}
  \cR_q = q^{\half\hJ_0 \otimes \hJ_0} \exp_{q^{-2}} \left(\left(1 -
q^{-2}\right)
  q^{\hJ_0/2} \hJ_+ \otimes q^{-\hJ_0/2} \hJ_-\right),
\end{equation}
where $\exp_q(\cX) = \sum_{n=0}^{\infty} \cX^n / \{n\}_q!$, $\{n\}_q! = \{n\}_q
\{n-1\}_q \cdots \{1\}_q$ if $n \in \N^+$, $\{0\}_q! = 1$,  and $\{n\}_q =
(1 - q^n) /
(1 - q)$. The generators of the standard U$_q$(sl(2)) algebra satisfy the
commutation relations
\begin{equation}
  q^{\hJ_0} \hJ_{\pm} q^{-\hJ_0} = q^{\pm2} \hJ_{\pm}, \qquad \left[\hJ_+, \hJ_-
  \right] = \left[\hJ_0\right]_q,
\end{equation}
where $[\cX]_q = (q^{\cal X} - q^{\cal X}) / (q - q^{-1})$.\par
%
%
{}Following Ref.~\cite{abdesselam98}, we introduce a transforming matrix~$M$,
singular in the $q \to 1$ limit, as
\begin{equation}
  M = \mbox{\rm E}_q\left(\eta \hJ_+\right),  \label{eq:M}
\end{equation}
where $\mbox{\rm E}_q(\cX) = \sum_{n=0}^{\infty} \cX^n / [n]_q!$, $[n]_q! =
[n]_q
[n-1]_q \cdots [1]_q$ if $n \in \N^+$, $[0]_q! = 1$,  and $\eta = h/(q-1)$.
The standard
$R_q^{j_1;j_2}$ matrix may now be subjected to a similarity transformation
followed by a $q \to 1$ limiting process:
\begin{equation}
  \tR_h^{j_1;j_2} = \lim_{q\to1} \left[\left(M^{-1}_{j_1} \otimes M^{-1}_{j_2}
  \right) R_q^{j_1;j_2} \left(M_{j_1} \otimes M_{j_2}\right)\right].
  \label{eq:lim-R}
\end{equation}
The matrix $\tR_h^{j_1;j_2}$ obtained by the above prescription is non-singular,
and may be seen to coincide, upon application of nonlinear
map~(\ref{eq:nl-map}),
with the corresponding $R_h^{j_1;j_2}$ obtained directly from the
expression~(\ref{eq:univ-R_h}) of the universal $\cR_h$ matrix.\par
%
%
The fundamental $j = \half$ representation of the U$_q$(sl(2)) algebra remains
undeformed, and the corresponding $\tR_h^{\half;j}$ matrix in~(\ref{eq:lim-R})
reads~\cite{abdesselam98}
\begin{equation}
  \tR_h^{\half;j} = \left(
      \begin{array}{cc}
      \tT & - {\ts\half} h \left(\tT + \tT^{-1}\right) J_0 + {\ts\half} h
           \left(\tT - \tT^{-1}\right) \\[0.2cm]
      0    & \tT^{-1}
      \end{array}\right).  \label{eq:tR-half}
\end{equation}
\par
%
%
Using similar techniques, the higher-dimensional $R_h$ matrices may also be
studied. The general features of the construction~(\ref{eq:lim-R}) may be
observed
in analysing the results for $(1 \otimes j)$, and $\left(\thalf \otimes
j\right)$
representations, respectively. For convenience, we choose the
$(2j+1)$-dimensional
irreducible representation of the U$_q$(sl(2)) algebra as follows:
\begin{eqnarray}
  \hJ_+ |jm\rangle & = & [j-m]_q [j+m+1]_q\, |j\; m+1\rangle, \qquad \hJ_-
          |jm\rangle = |j\; m-1\rangle, \nonumber \\
  \hJ_0 |jm\rangle & = & m |jm\rangle,  \label{eq:q-rep}
\end{eqnarray}
where $m=j$, $j-1$, $\ldots$, $-(j-1)$,~$-j$. The transforming matrix~$M$
in~(\ref{eq:M}) may now be obtained for $j=1$, and $j=\thalf$ representations as
\begin{equation}
  M_{j=1} = \left(
        \begin{array}{ccc}
        1 & [2]_q\, \eta & [2]_q\, \eta^2 \\[0.2cm]
        0 & 1              & [2]_q\, \eta \\[0.2cm]
        0 & 0              & 1
       \end{array}\right), \qquad
  M_{j=\thalf} = \left(
        \begin{array}{cccc}
        1 & [3]_q\, \eta & [3]_q!\, \eta^2 & [3]_q!\, \eta^3 \\[0.2cm]
        0 & 1              & [2]_q^2\, \eta     & [3]_q!\, \eta^2 \\[0.2cm]
        0 & 0              & 1                     & [3]_q\, \eta \\[0.2cm]

        0 & 0              & 0                     & 1
        \end{array}\right).  \label{eq:M-one}
\end{equation}
{}Following the recipe~(\ref{eq:lim-R}), the $\tR_h^{1;j}$ and
$\tR_h^{\thalf;j}$ matrices may be realized as
\begin{equation}
  \tR_h^{1;j} = \left(
        \begin{array}{ccc}
        \tT^2 & - h \left(\tT^2 + 1\right) J_0 & {\ts\half} h^2
\Bigl[\left(\tT +
               \tT^{-1}\right)^2 J_0^2 - 4 \left(\tT^2 - \tT^{-2}\right) \\
                 &                                                & +
4\left(\tT^{-2}
               - 1\right) J_0 \Bigr] \\[0.3cm]
        0       & 1                                             & - h
\left[\left(\tT^{-2} + 1\right) J_0
               + 2 \left(\tT^{-2} - 1\right)\right] \\[0.3cm]
        0       & 0                                             & \tT^{-2}
        \end{array}\right),  \label{eq:tR-one}
\end{equation}
and
\begin{equation}
  \tR_h^{\thalf;j} = \left(
        \begin{array}{cccc}
        \tT^3 & \tilde{A} & \tilde{B} & \tilde{C} \\[0.2cm]
        0       & \tT          & \tilde{D} & \tilde{E} \\[0.2cm]
        0       & 0             & \tT^{-1}  & \tilde{F} \\[0.2cm]
        0       & 0             & 0             & \tT^{-3}
        \end{array}\right),    \label{eq:tR-thalf}
\end{equation}
where
\begin{eqnarray}
  \tilde{A} & = & - {\ts\thalf} h \left[\tT^3 (J_0+1) + \tT (J_0-1)\right],
        \nonumber \\
  \tilde{B} & = & {\ts\thalf} h^2 \left[\tT^3 (J_0-1) (J_0+3) + 2 \tT (J_0-1)^2
        + \tT^{-1} (J_0+1)^2\right], \nonumber \\
  \tilde{C} & = & - {\ts\frac{3}{4}} h^3 \Bigl[ \tT^3 (J_0-3) (J_0+1) (J_0+5)
        + 3 \tT (J_0-3) (J_0-1) (J_0+1) \nonumber \\
  & & \mbox{} + 3 \tT^{-1} (J_0-3) (J_0+1)^2 + \tT^{-3} (J_0+1) (J_0+3) (J_0+5)
        \Bigr],\nonumber \\
  \tilde{D} & = & - 2h \left[\tT (J_0-1) + \tT^{-1} (J_0+1)\right], \nonumber \\
  \tilde{E} & = & {\ts\thalf} h^2 \left[\left(\tT + 2 \tT^{-1}\right)
(J_0-3) (J_0+1)
        + \tT^{-3} (J_0+3)^2\right],\nonumber \\
  \tilde{F} & = & - {\ts\thalf} h \left[\tT^{-1} (J_0-3) + \tT^{-3}
(J_0+3)\right].
 \end{eqnarray}
In the above computational procedure of the $\tR_h^{j_1;j_2}$ matrices, obtained
\`a la~(\ref{eq:lim-R}), the order of the singularity in the $q \to 1$
limit increases
by one for each step away from the diagonal, the rightmost corner element
being the
most singular. The singularities in the said $q \to 1$ limit, however, all
cancel
yielding finite results as in~(\ref{eq:tR-one}), and~(\ref{eq:tR-thalf}).\par
%
%
To prove the equivalence of the $\tR_h^{j_1;j_2}$ matrices obtained above as per
our contraction procedure, and the corresponding $R_h^{j_1;j_2}$ read directly
from the universal $\cR_h$ matrix~(\ref{eq:univ-R_h}), we first use the
map~(\ref{eq:nl-map}) to obtain a representation of the $h$-Jordanian
algebra~(\ref{eq:h-alg}). A $(2j+1)$-dimensional representation of the classical
algebra~(\ref{eq:class-alg})
\begin{eqnarray}
  J_+ |jm\rangle & = & (j-m) (j+m+1) |j\; m+1\rangle, \qquad J_- |jm\rangle =
          |j\; m-1\rangle, \nonumber \\
  J_0 |jm\rangle & = & m |j m\rangle,
\end{eqnarray}
now, via the nonlinear map~(\ref{eq:nl-map}), immediately furnishes the
corresponding $(2j+1)$-dimensional representation of the $h$-Jordanian
algebra~(\ref{eq:h-alg}). For the $j=\half$~case, the generators remain
undeformed.
For the $j=1$ and $j=\thalf$~cases, we list the representations of U$_h$(sl(2))
generators below:
\begin{eqnarray}
  \lefteqn{(j = 1)} \nonumber \\[0.2cm]
  & & X = \left(
         \begin{array}{ccc}
         0 & 2 & 0 \\[0.2cm]
         0 & 0 & 2 \\[0.2cm]
         0 & 0 & 0
         \end{array}\right), \qquad
  Y = \left(
         \begin{array}{ccc}
         0 & {\ts\half} h^2 & 0 \\[0.2cm]
         1 & 0                    & - {\ts\thalf} h^2 \\[0.2cm]
         0 & 1 & 0
         \end{array}\right), \nonumber \\[0.2cm]
  & & H = \left(
         \begin{array}{ccc}
         2 & 0 & -4 h^2 \\[0.2cm]
         0 & 0 & 0 \\[0.2cm]
         0 & 0 & -2
         \end{array}\right), \\[0.4cm]
  \lefteqn{(j = {\ts\thalf})} \nonumber \\[0.2cm]
  & & X = \left(
         \begin{array}{cccc}
         0 & 3 & 0 & -6 h^2 \\[0.2cm]
         0 & 0 & 4 & 0 \\[0.2cm]
         0 & 0 & 0 & 3 \\[0.2cm]
         0 & 0 & 0 & 0
         \end{array}\right), \qquad
  Y = \left(
         \begin{array}{cccc}
         0 & 0 & 0 & 0 \\[0.2cm]
         1 & 0 & 0 & 0 \\[0.2cm]
         0 & 1 & 0 & -6 h^2 \\[0.2cm]
         0 & 0 & 1 & 0
         \end{array}\right), \nonumber \\[0.2cm]
  & & H = \left(
         \begin{array}{cccc}
         3 & 0 & -6 h^2 & 0 \\[0.2cm]
         0 & 1 & 0        & -18 h^2 \\[0.2cm]
         0 & 0 & -1      & 0 \\[0.2cm]
         0 & 0 & 0        & -3
         \end{array}\right).
\end{eqnarray}
\par
%
%
Using the above representations in the expression~(\ref{eq:univ-R_h}) of the
universal $\cR_h$ matrix, we obtain
\begin{equation}
  R_h^{\half;j} = \left(
         \begin{array}{cc}
         T & -h H + {\ts\half} h \left(T - T^{-1}\right) \\[0.2cm]
         0 & T^{-1}
         \end{array}\right),  \label{eq:R-half}
\end{equation}
\begin{equation}
  R_h^{1;j} = \left(
         \begin{array}{ccc}
         T^2 & -2h TH & -2 h^2 \left[T^2 - T^{-2} + 2 TH \left(T^{-2} -
1\right) - (TH)^2
               T^{-2}\right] \\[0.2cm]
         0    & 1         & -2h \left(THT^{-2} - T^{-2} + 1\right) \\[0.2cm]
         0    & 0         & T^{-2}
         \end{array}\right),
\end{equation}
and
\begin{equation}
  R_h^{\thalf;j} = \left(
         \begin{array}{cccc}
         T^3 & A & B        & C \\[0.2cm]
         0    & T & D        & E \\[0.2cm]
         0    & 0 & T^{-1} & F \\[0.2cm]
         0    & 0 & 0        & T^{-3}
         \end{array}\right),  \label{eq:R-thalf}
\end{equation}
where
\begin{eqnarray}
  A & = & {\ts\thalf} h \left(T^2 - 1 - 2TH\right) T, \nonumber \\
  B & = & - {\ts\thalf} h^2 \left[3T^3 - 2T - T^{-1} - 4TH \left(T^2 - 1 +
TH\right)
         T^{-1}\right], \nonumber \\
  C & = & - 3h^3 \Bigl[{\ts\frac{15}{4}} \left(T^3 - T^{-3}\right) - {\ts\frac{9}{4}}
         \left(T - T^{-1}\right) - {\ts\frac{9}{2}} THT - 9THT^{-1} \nonumber \\
  & & \mbox{} + {\ts\frac{23}{2}} THT^{-3} - 9 (TH)^2 \left(T^{-3} -
T^{-1}\right)
         + 2 (TH)^3 T^{-3}\Bigr], \nonumber \\
  D & = & - 2h \left(T^2 - 1 + 2TH\right) T^{-1}, \nonumber \\
  E & = & - {\ts\thalf} h^2 \left\{3 \left(T + 2T^{-1} - 3T^{-3}\right) -
4TH \left[
         3 \left(T^2 - 1\right) + TH\right] T^{-3}\right\}, \nonumber \\
  F & = & - {\ts\thalf} h \left[3 \left(T^2 - 1\right) + 2TH\right] T^{-3}.
\end{eqnarray}
Upon applying the map~(\ref{eq:nl-map}), the matrices~(\ref{eq:tR-half}),
(\ref{eq:tR-one}), and~(\ref{eq:tR-thalf}), obtained using our contraction
procedure,
may be shown to be equal to the corresponding
matrices~(\ref{eq:R-half})--(\ref{eq:R-thalf}), obtained directly using the
expression~(\ref{eq:univ-R_h}) of the universal $\cR_h$ matrix. This process may
be continued indefinitely for arbitrary representations signifying that the
arbitrary
$h$-Jordanian matrices $R_h^{j_1;j_2}$ may be recovered, via the contraction
route~(\ref{eq:lim-R}), from the corresponding standard $R_q^{j_1;j_2}$
matrices.\par
%
%
\section{On U$_{h,\alpha}$(gl(2)) and Its $R_{h,\alpha}$ Matrices}
\label{sec:R_{h,alpha}}
\setcounter{equation}{0}
We now consider the two-parametric $(h, \alpha)$-Jordanian deformed algebra
U$_{h,\alpha}$(gl(2))${}\sim{}$U$_{h,\alpha}(\mbox{\rm sl(2)} \oplus \mbox{\rm
u(1)})$, generated by $(X, Y, H, Z)$~\cite{agha93, aneva}. The
generator~$Z$ of the
u(1) algebra is a primitive central element. Our choice of the deformation
parameters is a little different from others. Our results might be rewritten in
terms of the customary choice of deformation parameters $h_{\pm} = h (1 \pm
\alpha)$. By studying the duality relation between the Hopf algebras
Fun$_{h,\alpha}$(GL(2)) and U$_{h,\alpha}$(gl(2)), the Hopf structure of
the latter
has been investigated by Aneva {\sl et al.}~\cite{aneva}. After a suitable
choice of
parameters, these authors conclude that the role of the second parameter may be
confined to the coalgebra and the antipode alone.\par
%
%
Here we use the general formalism developed by Reshetikhin~\cite{reshe} for
introducing multiple deformation parameters. This allows us to obtain the
universal $\cR_{h,\alpha}$ matrix of the two-parametric U$_{h,\alpha}$(gl(2))
algebra. To this end, we construct a twist operator
\begin{equation}
  F_{\alpha} = \exp [\alpha h (X \otimes Z - Z \otimes X)].
\end{equation}
{}Following Ref.~\cite{reshe}, the coproducts of U$_h$(gl(2)) and
U$_{h,\alpha}$(gl(2)) may be related by
\begin{equation}
  \Delta_{h,\alpha}(\cX) = F_{\alpha} \Delta_h(\cX) F_{\alpha}^{-1} \qquad
\forall
  \cX \in \{X, Y, H, Z\}.
\end{equation}
Explicitly, we obtain
\begin{eqnarray}
  \Delta_{h,\alpha}(X) & = & X \otimes 1 + 1 \otimes X, \nonumber \\
  \Delta_{h,\alpha}(Y) & = & Y \otimes T + T^{-1} \otimes Y + \alpha h \left(H
          \otimes ZT - ZT^{-1} \otimes H\right) \nonumber\\
  & & \mbox{} - {\ts\half} \alpha^2 h \left[Z^2 T^{-1} \otimes \left(T -
T^{-1}\right)
          + \left(T - T^{-1}\right) \otimes Z^2 T\right], \nonumber \\
  \Delta_{h,\alpha}(H) & = & H \otimes T + T^{-1} \otimes H + \alpha
\left[Z T^{-1}
          \otimes \left(T - T^{-1}\right) - \left(T - T^{-1}\right) \otimes
Z T\right],
          \nonumber \\
  \Delta_{h,\alpha}(Z) & = & Z \otimes 1 + 1 \otimes Z.
\end{eqnarray}
The antipode and the counit maps, however, remain unaltered. The authors of
Ref.~\cite{aneva} obtained the antipode maps for U$_{h,\alpha}$(gl(2))
depending on
both the deformation parameters. This result is, however, an artifact of their
choice of the corresponding coproduct maps in a asymmetric way.\par
%
%
The universal $\cR_{h,\alpha}$ matrix of the U$_{h,\alpha}$(gl(2)) algebra
may also
be determined following the Reshetikhin procedure~\cite{reshe}:
\begin{equation}
  \cR_{h,\alpha} = F_{\alpha}^{-1} \cR_h F_{\alpha}^{-1},
  \label{eq:univ-R_halpha}
\end{equation}
which, by construction, satisfies the properties~(\ref{eq:R-prop}), and the
triangularity condition. Explicitly, $\cR_{h,\alpha}$ reads
\begin{eqnarray}
  \cR_{h,\alpha} & = & \exp[- 2h\alpha (X \otimes Z - Z \otimes X)]
\exp\left[- h
           \left(X \otimes TH + \alpha ZX \otimes
\left(T^2-1\right)\right)\right]
           \nonumber \\
  & & \mbox{} \times \exp\left[h \left(TH \otimes X - \alpha \left(T^2-1\right)
           \otimes ZX\right)\right].
\end{eqnarray}
\par
%
%
The above universal $\cR_{h,\alpha}$ matrix of the $(h,\alpha)$-Jordanian
algebra
U$_{h,\alpha}$(gl(2)) may also be obtained by the contraction technique
discussed
in Sect.~\ref{sec:R_h}. Our starting point here is the universal
$\cR_{q,\lambda}$
matrix~\cite{chakra} of the two-parametric $(q,\lambda)$-deformed standard
U$_{q,\lambda}$(gl(2)) algebra~\cite{schirrmacher}. Following Reshetikhin
procedure~\cite{reshe}, it was observed~\cite{chakra} that the universal
$\cR_{q,\lambda}$ matrix of the U$_{q,\lambda}$(gl(2)) algebra may be related to
the universal $\cR_q$ matrix of the standard $q$-deformed U$_q$(sl(2))
algebra via
a twist operator:
\begin{equation}
  \cR_{q,\lambda} = \hF_{\lambda}^{-1} \cR_q \hF_{\lambda}^{-1},
  \label{eq:univ-R_qlambda}
\end{equation}
where $\hF_{\lambda}$ is given by
\begin{equation}
  \hF_{\lambda} = \lambda^{\half\left(\hJ_0 \otimes \hZ - \hZ \otimes \hJ_0
  \right)} = q^{\frac{\alpha}{2}\left(\hJ_0 \otimes \hZ - \hZ \otimes \hJ_0
  \right)}.  \label{eq:hF}
\end{equation}
In the second equality in~(\ref{eq:hF}), we assume without any loss of
generality
$\lambda = q^{\alpha}$.\par
%
%
After taking into account the discussions in Sect.~\ref{sec:R_h} about the
convertibility, via the contraction process, of the standard $\cR_q$ matrix
of the
U$_q$(sl(2)) algebra into the Jordanian $\cR_h$ matrix of the U$_h$(sl(2))
algebra,
a comparison between~(\ref{eq:univ-R_halpha}) and~(\ref{eq:univ-R_qlambda})
makes it obvious that we need to show
\begin{equation}
  F_{\alpha}^{-1} = \lim_{q\to1} \left(M^{-1} \otimes M^{-1}\right)
  \hF_{\lambda}^{-1} (M \otimes M).  \label{eq:F-hF}
\end{equation}
Equation~(\ref{eq:F-hF}) may be readily proved as an operator identity.
Using~(\ref{eq:hF}), and following the method developed in
Ref.~\cite{abdesselam98}, we may show that the rhs of~(\ref{eq:F-hF}) assumes
the form
\begin{eqnarray}
  \mbox{rhs} & = & \tT^{-\alpha z_2} \otimes \tT^{\alpha z_1} \nonumber \\
  & = & \exp[- \alpha h (X \otimes Z - Z \otimes X)] = \mbox{lhs}.
\label{eq:rhs-lhs}
\end{eqnarray}
In the first equality in~(\ref{eq:rhs-lhs}), the central generator~$\hZ$ of the
U$_{q,\lambda}$(gl(2)) algebra is assumed to have values~$z_1$ and~$z_2$ in the
two sectors of the tensor product space. In the second equality
in~(\ref{eq:rhs-lhs}), we have used the map~(\ref{eq:nl-map}), and the
identity map
$Z = \hZ$. This completes our demonstration that the $(h, \alpha)$-Jordanian
$R_{h,\alpha}^{j_1,z_1;j_2,z_2}$ matrix of arbitrary coloured
representation $(j_1,
z_1 \otimes j_2, z_2)$ may be recovered through the present contraction process
from the corresponding  two-parametric standard
$R_{q,\lambda}^{j_1,z_1;j_2,z_2}$ matrix.\par
%
%
As an application of our method, we here obtain the matrices
$R_{h,\alpha}^{\half,z_1;1,z_2}$ and $R_{h,\alpha}^{1,z_1;\half,z_2}$, which
satisfy appropriate coloured Yang-Baxter equation~(\ref{eq:YB}) with the central
generator assuming different values in different sectors of the tensor product
space:
\begin{equation}
  R_{h,\alpha}^{\half,z_1;1,z_2} = \left(
         \begin{array}{cc}
         A' & B' \\[0.2cm]
         0  & C'
         \end{array}\right), \qquad
  R_{h,\alpha}^{1,z_1;\half,z_2} = \left(
         \begin{array}{ccc}
         A'' & B'' & C'' \\[0.2cm]
         0  &  D'' & E'' \\[0.2cm]
         0  &  0  & F''
         \end{array}\right),  \label{eq:col-R}
\end{equation}
where $A'$, $B'$, $C'$ are $3 \times 3$ matrices,
\begin{eqnarray}
  A' & = & \left(
        \begin{array}{ccc}
        1 & 2h (1 + 2\alpha z_1) & 2h^2 (1 + 2\alpha z_1)^2 \\[0.2cm]
        0 & 1                               & 2h (1 + 2\alpha z_1) \\[0.2cm]
        0 & 0                               & 1
        \end{array}\right), \nonumber \\[0.2cm]
  B' & = & \left(
        \begin{array}{ccc}
        - 2h (1 + \alpha z_2) & 2h^2 (1 - 2\alpha z_1  & 4h^3 (1 +
              2\alpha z_1 - \alpha z_2 \\
                                         &  - 4\alpha^2 z_1 z_2)   & -
4\alpha^3 z_1^2 z_2)
              \\[0.2cm]
        0                               & -2h\alpha z_2 & 2h^2 (1
              + 2\alpha z_1 - 4\alpha^2 z_1z_2) \\[0.2cm]
        0                               & 0 & 2h (1 -
              \alpha z_2)
        \end{array}\right), \nonumber \\[0.2cm]
  C' & = & \left(
        \begin{array}{ccc}
        1 & - 2h (1 - 2\alpha z_1) & 2h^2 (1 - 2\alpha z_1)^2 \\[0.2cm]
        0 & 1                                  & - 2h (1 - 2\alpha z_1)
\\[0.2cm]
        0 & 0                                  & 1
        \end{array}\right),
\end{eqnarray}
while $A''$, $B''$, $C''$, $D''$, $E''$, $F''$ are $2 \times 2$ matrices,
\begin{eqnarray}
  A'' & = & \left(
         \begin{array}{cc}
         1 & 2h (1 + \alpha z_1) \\[0.2cm]
         0 & 1
         \end{array}\right), \nonumber \\[0.2cm]
  B'' & = & \left(
         \begin{array}{cc}
         - 2h (1 + 2\alpha z_2) & 2h^2 (1 - 2\alpha z_2 - 4\alpha^2 z_1
z_2) \\[0.2cm]
         0                                 & 2h (1 - 2\alpha z_2)
         \end{array}\right), \nonumber \\[0.2cm]
  C'' & = & \left(
         \begin{array}{cc}
         2h^2 (1 + 2\alpha z_2)^2 & - 4h^3 (1 + 2\alpha z_2 - \alpha z_1 -
4\alpha^3
               z_1 z_2^2) \\[0.2cm]
         0 & 2h^2 (1 - 2\alpha z_2)^2
         \end{array}\right), \nonumber \\[0.2cm]
  D'' & = & \left(
         \begin{array}{cc}
         1 & 2h\alpha z_1 \\[0.2cm]
         0 & 1
         \end{array}\right), \nonumber \\[0.2cm]
   E'' & = & \left(
         \begin{array}{cc}
         - 2h (1 + 2\alpha z_2) & 2h^2 (1 + 2\alpha z_2 - 4\alpha^2 z_1
z_2) \\[0.2cm]
         0                                  & 2h (1 - 2\alpha z_2)
         \end{array}\right), \nonumber \\[0.2cm]
  F'' & = & \left(
         \begin{array}{cc}
         1 & - 2h (1 - \alpha z_1) \\[0.2cm]
         0 & 1
         \end{array}\right).
\end{eqnarray}
Comparing the two matrices in~(\ref{eq:col-R}), it is evident that an  `exchange
symmetry' between the two sectors of the tensor product space holds:
\begin{equation}
  \left(R_{h,\alpha}^{j_1,z_1;j_2,z_2}\right)_{km,ln} =
  \left(R_{-h,\alpha}^{j_2,z_2;j_1,z_1}\right)_{mk,nl}.
\end{equation}
\par
%
%
\section{The Coloured $T_{h,\alpha}^{j,z}$ Matrices of the
Fun$_{h,\alpha}$(GL(2))
Algebra}
\label{sec:T_{h,alpha}}
\setcounter{equation}{0}
In this section, we turn our attention to the determination of the
$(h,\alpha)$-Jordanian $T$ matrices in arbitrary coloured representation. The
contraction method described earlier may also be employed here. The
technique may
be summarized as follows. Starting with a two-parametric $(q,\lambda)$-standard
$T_{q,\lambda}^{j,z}$ matrix of the Fun$_{q,\lambda}$(GL(2)) algebra in an
arbitrary representation given by~$j$, and the colour parameter~$z$, we
perform a
suitable similarity transformation by the matrix~$M_j$ defined in~(\ref{eq:M}).
Even though the transforming matrix is singular in the $q\to1$ limit, the
singularities all systematically cancel for the transformed matrix
\begin{equation}
  \tT^{j,z} = M_j^{-1} T_{q,\lambda (=q^{\alpha})}^{j,z} M_j,  \label{eq:tT}
\end{equation}
which yields finite elements of the corresponding two-parametric
$(h,\alpha)$-Jordanian $T_{h,\alpha}^{j,z}$ matrix for the $(j,z)$ representation:
\begin{equation}
  T_{h,\alpha}^{j,z} = \lim_{q\to1} \tT^{j,z}.  \label{eq:T}
\end{equation}
The $T_{h,\alpha}^{j,z}$ matrix in~(\ref{eq:T}) satisfies the coloured inverse
scattering equation~(\ref{eq:RTT}). For the fundamental representation $(j
= {\ts
\half}, z = {\ts \half})$, the method was previously employed~\cite{agha95}. We
first reproduce these results here for the notational purpose.\par
%
%
The generating elements $(\ha, \hb, \hc, \hd)$ of the two-parametric standard
Fun$_{q,\lambda}$(GL(2)) algebra, defined as $\hT_{q,\lambda}^{\half,\half} =
\left(\begin{array}{cc} \ha & \hb \\ \hc & \hd \end{array}\right)$, satisfy the
following commutation relations
\begin{eqnarray}
  \ha \hb & = & q^{-1} \lambda^{-1} \hb \ha, \qquad \ha \hc = q^{-1}
\lambda \hc \ha,
          \qquad \hb \hd = q^{-1} \lambda \hd \hb, \qquad \hc \hd = q^{-1}
\lambda^{-1}
          \hd \hc, \nonumber \\
  \left[\ha, \hd\right] & = & \left(q^{-1} - q\right) \lambda^{-1} \hb \hc,
\qquad
          \lambda^{-1} \hb \hc = \lambda \hc \hb,  \label{eq:Funq-alg}
\end{eqnarray}
and the matrix coproduct rule
$\Delta\left(\hT_{q,\lambda}^{\half,\half}\right) =
\hT_{q,\lambda}^{\half,\half} \dotimes \hT_{q,\lambda}^{\half,\half}$. The
determinant $\hD = \ha \hd - q^{-1} \lambda^{-1} \hb \hc$ has the algebraic
properties
\begin{equation}
  \left[\ha, \hD\right] = 0, \qquad \hb \hD = \lambda^2 \hD \hb, \qquad \hc
\hD =
  \lambda^{-2} \hD \hc, \qquad \left[\hd, \hD\right] = 0,
\end{equation}
and is endowed with a group-like coproduct $\Delta\left(\hD\right) = \hD \otimes
\hD$.\par
%
%
The elements of the transformed $\tT^{\half,\half} =
\left(\begin{array}{cc} \ta &
\tb \\ \tc & \td \end{array}\right)$ matrix, obtained
following~(\ref{eq:tT}), are
related to the old basis set as
\begin{equation}
  \ha = \ta + \eta \tc, \qquad \hb = \tb - \eta \left(\ta - \td\right) -
\eta^2 \tc,
  \qquad \hc = \tc, \qquad \hd = \td - \eta \tc.  \label{eq:ha-ta}
\end{equation}
Staying away from the singularity in~(\ref{eq:ha-ta}) in the $q\to1$ limit, we
recast the algebra~(\ref{eq:Funq-alg}) as
\begin{eqnarray}
  \ta \tb - q^{-1-\alpha} \tb \ta & = & \eta \left(1 -
q^{-1-\alpha}\right)\left(\ta^2
          - \tD\right), \nonumber \\
  \ta \tc - q^{-1+\alpha} \tc \ta & = & - \eta \left(1 -
q^{-1+\alpha}\right) \tc^2,
          \nonumber \\
  \tb \td - q^{-1+\alpha} \td \tb & = & - \eta \left(1 - q^{-1+\alpha}\right)
          \left(\td^2 - \tD\right), \nonumber \\
  \tc \td - q^{-1-\alpha} \td \tc & = & \eta \left(1 - q^{-1-\alpha}\right)
\tc^2,
          \nonumber \\
  \left[\ta, \td\right] & = & \eta \left(1 - q^{-1-\alpha}\right) \ta \tc - \eta
          \left(q^{1-\alpha} - 1\right) \td \tc + \left(q^{-1} - q\right)
q^{-\alpha} \tb
          \tc, \nonumber \\
  q^{-\alpha} \tb \tc - q^{\alpha} \tc \tb & = & - \eta \left(q -
q^{-\alpha}\right)
          \ta \tc - \eta \left(q - q^{\alpha}\right) \tc \td,
\label{eq:Funq-talg}
\end{eqnarray}
where the determinant $\tD = \ta \td - q^{-1-\alpha} \tb \tc - \eta \left(1 -
q^{-1-\alpha}\right) \ta \tc$ satisfies the commutation relations
\begin{eqnarray}
  \left[\ta, \tD\right] & = & \eta \left(q^{2\alpha} - 1\right) \tc \tD, \qquad
          \tb \tD - q^{2\alpha} \tD \tb = \eta \left(q^{2\alpha} - 1\right)
\left(\tD \td
          - \ta \tD\right), \nonumber \\
  \tc \tD - q^{-2\alpha} \tD \tc & = & 0, \qquad   \left[\td, \tD\right] =
- \eta
          \left(q^{2\alpha} - 1\right) \tc \tD.  \label{eq:Funq-talgbis}
\end{eqnarray}
The elements obey the matrix coproduct rule
$\Delta\left(\tT^{\half,\half}\right) = \tT^{\half,\half} \dotimes
\tT^{\half,\half}$.\par
%
%
Defining the generating elements $(a, b, c, d)$ of the Fun$_{h,\alpha}$(GL(2))
algebra as $T_{h,\alpha}^{\half,\half} = \left(\begin{array}{cc} a & b \\ c & d
\end{array}\right)$, we now pass to the $q\to1$ limit in~(\ref{eq:Funq-talg})
and~(\ref{eq:Funq-talgbis}), and obtain the following:
\begin{eqnarray}
  [a, b] & = & h (1 + \alpha) \left(a^2 - D\right), \qquad [a, c] = - h (1
- \alpha) c^2,
          \nonumber \\[0cm]
  [b, d] & = & - h (1 - \alpha) \left(d^2 - D\right), \qquad [c, d] = h (1
+ \alpha) c^2,
          \nonumber \\[0cm]
  [a, d] & = & h (1 + \alpha) ac - h (1 - \alpha) dc, \qquad [b, c] = - h
(1 + \alpha) ac
          - h (1 - \alpha) cd,  \label{eq:Funh-alg}
\end{eqnarray}
where the determinant $D = ad - bc - h(1+ \alpha) ac$ satisfies the commutation
relations
\begin{equation}
  [a, D] = 2h \alpha cD, \quad [b, D] = 2h \alpha (Dd - aD), \quad [c, D] =
0, \quad
  [d, D] = - 2h \alpha cD.
\end{equation}
The coproduct rules read
\begin{equation}
  \Delta \left(\begin{array}{cc}
                      a & b \\[0.2cm]
                      c & d
                      \end{array}\right)
  = \left(\begin{array}{cc}
               a & b \\[0.2cm]
               c & d
               \end{array}\right) \dotimes
    \left(\begin{array}{cc}
               a & b \\[0.2cm]
               c & d
               \end{array}\right).  \label{eq:Funh-coalg}
\end{equation}
It is interesting to observe that the algebraic relations~(\ref{eq:Funh-alg})
exhibit an automorphism
\begin{equation}
  \rho(a) = d, \qquad \rho(b) = b, \qquad \rho(c) = c, \qquad \rho(d) = a,
\end{equation}
where the deformation parameters are also redefined according to $h \to h$,
$\alpha \to - \alpha$. The map~$\rho$ is also an anti-automorphism of the
coalgebra rules~(\ref{eq:Funh-coalg}):
\begin{equation}
  \sigma \circ (\rho \otimes \rho) \circ \Delta(\rho(\cX)) = \Delta(\cX) \qquad
  \forall \cX \in \{a, b, c, d\}.
\end{equation}
\par
%
%
As the first example of our proposed generalization --- as indicated
in~(\ref{eq:tT}) and~(\ref{eq:T}) --- of the above construction method, we now
consider $(h, \alpha)$-Jordanian $T_{h,\alpha}^{j=\half,z}$ matrix for an
arbitrary
colour parameter~$z$. The corresponding $(q, \lambda)$-standard
$T_{q,\lambda}^{\half,z}$ matrix for $\left(j = {\ts\half}, z\right)$
representation
may be obtained~\cite{jaga} from the known expression~\cite{fronsdal} of the
universal $\cT_{q,\lambda}$ matrix acting as a dual form between the
Fun$_{q,\lambda}$(GL(2)) and U$_{q,\lambda}$(gl(2)) algebras,
\begin{equation}
  \cT_{q,\lambda} = \exp_{q^2}\left(\hgamma \hcJ_-\right) \exp\left[\halpha
\left(
  \frac{\hJ_0}{2} + \hZ\right) + \hdelta \left(\frac{\hJ_0}{2} -
\hZ\right)\right]
  \exp_{q^{-2}}\left(\hbeta \hcJ_+\right),  \label{eq:univ-T}
\end{equation}
where $\hcJ_+ = \hJ_+ q^{(\hJ_0 + 1)/2} \lambda^{\hZ - \half}$ and $\hcJ_- =
q^{-(\hJ_0 + 1)/2} \lambda^{\hZ - \half} \hJ_-$. The redefined generating
elements
of the Fun$_{q,\lambda}$(GL(2)) algebra are given by
\begin{equation}
  \ha = e^{\halpha}, \qquad \hb = e^{\halpha} \hbeta, \qquad \hc = \hgamma
  e^{\halpha}, \qquad \hd = \hgamma e^{\halpha} \hbeta + e^{-\hdelta}.
\end{equation}
{}For the $\left(j = {\ts\half}, z\right)$ representation, the $(q,
\lambda)$-deformed standard $\hT_{q,\lambda}^{\half,z}$ matrix reads~\cite{jaga}
\begin{equation}
  \hT_{q,\lambda}^{\half,z} = \hD^{z-\half} \left(
      \begin{array}{cc}
      \ha                                  & \lambda^{z-\half} \hb \\[0.2cm]
      \lambda^{-(z-\half)} \hc & \hd
      \end{array}\right).
\end{equation}
\par
%
%
Using the transformation~(\ref{eq:tT}) and the redefinition of the
variables~(\ref{eq:ha-ta}), we now obtain
\begin{equation}
  \tT^{\half,z} = M_{\half}^{-1} \hT_{q,\lambda(= q^{\alpha})}^{\half,z}
M_{\half},
\end{equation}
as
\begin{equation}
  \tT^{\half,z} = \hD^{z-\half} \left(
        \begin{array}{cc}
        \ta + \eta \left(1 - q^{-\alpha(z-\half)}\right) \tc & q^{\alpha
(z-\half)} \hb
            - \eta \left(q^{\alpha (z-\half)} - 1\right) \left(\ta -
\td\right) \\
                                                                                        & - \eta^2 \left(q^{\alpha
            (z-\half)} + q^{-\alpha (z-\half)} - 2\right) \tc \\[0.2cm]
        q^{-\alpha (z-\half)} \tc
& \td - \eta \left(1 -
            q^{-\alpha (z-\half)}\right) \tc
        \end{array}\right).  \label{eq:tT-half}
\end{equation}
The limiting process $q\to1$ in~(\ref{eq:tT-half}) now yields the desired
$(h,\alpha)$-Jordanian $T_{h,\alpha}^{\half,z}$ matrix for the coloured
$\left(j = {\ts\half}, z\right)$ representation,
\begin{equation}
  T_{h,\alpha}^{\half,z} = D^{z-\half} \left(
        \begin{array}{cc}
        a + h \alpha (z-{\ts\half}) c & b - h \alpha (z-{\ts\half}) (a-d) -
h^2 \alpha^2
             (z-{\ts\half})^2 c \\[0.2cm]
        c                                           & d - h \alpha
(z-{\ts\half}) c
        \end{array}\right),  \label{eq:T-half}
\end{equation}
which satisfies the appropriate coloured inverse scattering
equation~(\ref{eq:RTT}). Substitution of the $z = \half$ value of the colour
parameter obviously reduces~(\ref{eq:T-half}) to the fundamental
representation.\par
%
%
As another demonstration of our procedure, we explicitly obtain the
$(h,\alpha)$-Jordanian $T_{h,\alpha}^{1,z}$ matrix for $j=1$ representation
at an
arbitrary value of the colour parameter~$z$. As previously, our starting
point here
is the $(q,\lambda)$-standard $\hT_{q,\lambda}^{1,z}$ matrix of the $(1,z)$
representation. Using the representation~(\ref{eq:q-rep}) and the universal
$\cT_{q,\lambda}$ matrix~(\ref{eq:univ-T}), we obtain
\begin{equation}
  \hT_{q,\lambda}^{1,z} = \hD^{z-1} \left(
        \begin{array}{ccc}
        \ha^2                                                 & [2]_q\,
q^{\half} \lambda^{z-\half} \ha
              \hb & [2]_q\, \lambda^{2(z-1)} \hb^2 \\[0.2cm]
        q^{\half} \lambda^{-z+\half} \ha \hc & \ha \hd + q \lambda^{-1} \hb \hc
                    & [2]_q\, q^{\half} \lambda^{z-\thalf} \hb \hd \\[0.2cm]
        [2]_q^{-1} \lambda^{-2(z-1)} \hc^2    & q^{\half}
\lambda^{-z+\thalf} \hc \hd
                    & \hd^2
        \end{array}\right).
\end{equation}
\par
%
%
Using the construction of the transforming matrix~$M_{j=1}$ in~(\ref{eq:M-one}),
we first perform the similarity transformation~(\ref{eq:tT}), and
using~(\ref{eq:ha-ta}) reexpress the elements of the relevant $\tT^{1,z}$
matrix in
terms of the variables $(\ta, \tb, \tc, \td)$. The transformed matrix is
now free
from singularities at all orders as~$q\to1$. Consequently, in the said
$q\to1$~limit, we now explicitly obtain finite elements of the corresponding
$(h,\alpha)$-Jordanian $T_{h,\alpha}^{1,z}$ matrix for the $(j=1, z)$
representation,
\begin{equation}
  T_{h,\alpha}^{1,z} = D^{z-1} \left(
         \begin{array}{ccc}
         U^z_{11} & U^z_{12} &U^z_{13} \\[0.2cm]
         U^z_{21} & U^z_{22} &U^z_{23} \\[0.2cm]
         U^z_{31} & U^z_{32} &U^z_{33}
         \end{array}\right),  \label{eq:T-one}
\end{equation}
where
\begin{eqnarray}
  U^z_{11} & = & a^2 + 2h\alpha (z-1) ac + {\ts \frac{1}{4}} h^2 [1 + \alpha
        (2z-3)]^2 c^2, \nonumber \\
  U^z_{12} & = & 2ab - h [1 + \alpha (2z-1)] \left(a^2 - ad\right) - h [1 -
\alpha
        (2z-3)] bc \nonumber \\
  & & \mbox{} - h^2 \left[1 + 2\alpha + \alpha^2 \left(4z^2 - 8z +
5\right)\right]
        ac + {\ts\half} h^2 [1 + \alpha (2z-3)]^2 cd \nonumber \\
  & & \mbox{} - {\ts \frac{1}{4}} h^3 \Bigl[1 + \alpha (6z-7) + \alpha^2
\left(12z^2 -
        28z + 15\right) \nonumber \\
  & & \mbox{} + \alpha^3 \left(8z^3 - 28z^2 + 30z -9\right)\Bigr] c^2,
        \nonumber \\
  U^z_{13} & = & 2b^2 - 4h\alpha (z-1) (ab-bd) - {\ts\half} h^2 \left[3 -
2\alpha
        (2z-1) - \alpha^2 (2z-1)^2\right] a^2 \nonumber \\
  & & \mbox{} + h^2 \left[1 - 4\alpha (z-1) - \alpha^2 \left(4z^2 - 8z +5\right)
        \right] ad \nonumber \\
  & & \mbox{} - h^2 \left[1 - 4\alpha (z-1) + \alpha^2 \left(4z^2 - 8z +3\right)
        \right] bc + {\ts\half} h^2 [1 + \alpha (2z-3)]^2 d^2 \nonumber \\
  & & \mbox{} - h^3 \left[1 - \alpha (z-2) - \alpha^2 (6z-5) - \alpha^3
\left(4z^3
        - 12z^2 + 13z - 4\right)\right] ac \nonumber \\
  & & \mbox{} - h^3 \left[\alpha (z-1) + 2\alpha^2 (2z^2 - 5z + 3) + \alpha^3
        \left(4z^3 - 16z^2 + 21z - 9\right)\right] cd \nonumber \\
  & & \mbox{} - {\ts \frac{1}{8}} h^4 \Bigl[3 + 8\alpha (z-2) - 2\alpha^2
\left(4z^2
        - 7\right) - 8\alpha^3 \left(4z^3 - 12z^2 + 11z - 3\right) \nonumber \\
  & & \mbox{} - \alpha^4 \left(16z^4 - 64z^3 + 88z^2 - 48z + 9\right)\Bigr] c^2,
        \nonumber \\
  U^z_{21} & = & ac + {\ts\half} h [1 + \alpha (2z-3)] c^2, \nonumber \\
  U^z_{22} & = & ad + bc - 2h\alpha (z-1) ac + h [1 + \alpha (2z-3)] cd
\nonumber \\
  & & \mbox{} - {\ts\half} h^2 \left[1 + 4\alpha (z-1) + \alpha^2
\left(4z^2 - 8z + 3
        \right)\right] c^2, \nonumber \\
  U^z_{23} & = & 2bd + h [1 + \alpha (2z-3)] \left(d^2 - ad\right) + h [1 -
\alpha
        (2z-1)] bc \nonumber \\
  & & \mbox{} - {\ts\half} h^2 \left[1 - 2\alpha - \alpha^2 \left(4z^2 - 8z
+3\right)
        \right] ac - 2h^2 \alpha \Bigl[z - 1 \nonumber \\
  & & \mbox{} + \alpha \left(2z^2 - 5z + 3\right)\Bigr] cd - {\ts
\frac{1}{4}} h^3
        \Bigl[3 + \alpha (2z-7) \nonumber \\
  & & \mbox{} - \alpha^2 \left(12z^2 - 20z + 7\right) - \alpha^3 \left(8z^3
- 20z^2 +
        14z - 3\right)\Bigr] c^2, \nonumber \\
  U^z_{31} & = & {\ts\half} c^2, \nonumber \\
  U^z_{32} & = & cd - {\ts\half} h [1 + \alpha (2z-1)] c^2, \nonumber \\
  U^z_{33} & = & d^2 - 2h\alpha (z-1) cd - {\ts \frac{1}{4}} h^2 \left[3 -
2\alpha
        (2z-1) - \alpha^2 (2z-1)^2\right] c^2.  \label{eq:T-one-bis}
\end{eqnarray}
\par
%
%
{}For the usual choice of the colour parameter~$z = \half$, the
$T_{h,\alpha}^{j=1, z=\half}$ matrix is obtained from~(\ref{eq:T-one})
and~(\ref{eq:T-one-bis}):
\begin{equation}
  T_{h,\alpha}^{1,\half} = D^{-1/2} \left(
         \begin{array}{ccc}
         U_{11} & U_{12} &U_{13} \\[0.2cm]
         U_{21} & U_{22} &U_{23} \\[0.2cm]
         U_{31} & U_{32} &U_{33}
         \end{array}\right),  \label{eq:T-one-half}
\end{equation}
where
\begin{eqnarray}
  U_{11} & = & a^2 - h\alpha ac + {\ts \frac{1}{4}} h^2 (1 - 2\alpha)^2
c^2, \nonumber
        \\
  U_{12} & = & 2ab - h \left(a^2 - ad\right) - h (1 + 2\alpha) bc - h^2
\left(1 +
        2\alpha + 2\alpha^2 \right) ac \nonumber \\
  & & \mbox{} + {\ts\half} h^2 (1 - 2\alpha )^2 cd - {\ts \frac{1}{4}} h^3
        (1 - 2\alpha)^2 c^2 \nonumber \\
  U_{13} & = & 2b^2 + 2h\alpha (ab-bd) - {\ts\thalf} h^2 a^2 + h^2 \left(1
+ 2\alpha
        - 2\alpha^2\right) ad \nonumber \\
  & & \mbox{} - h^2 (1 + 2\alpha) bc + {\ts\half} h^2 (1 - 2\alpha)^2 d^2
        - {\ts\half} h^3 \left(2 + 3\alpha + 4\alpha^2\right) ac \nonumber \\
  & & \mbox{} + {\ts\half} h^3 \alpha (1 - 2\alpha)^2 cd  - {\ts
\frac{3}{8}} h^4
        (1 - 2\alpha)^2 c^2, \nonumber \\
  U_{21} & = & ac + {\ts\half} h (1 - 2\alpha) c^2, \nonumber \\
  U_{22} & = & ad + bc + h\alpha ac + h (1 - 2\alpha) cd - {\ts\half} h^2
(1 - 2\alpha)
        c^2, \nonumber \\
  U_{23} & = & 2bd + h (1 - 2\alpha) \left(d^2 - ad\right) + h bc -
{\ts\half} h^2
        (1 - 2\alpha) ac \nonumber \\
  & & \mbox{} + h^2 \alpha (1 - 2\alpha) cd - {\ts \frac{3}{4}} h^3 (1 -
2\alpha) c^2,
        \nonumber \\
  U_{31} & = & {\ts\half} c^2, \nonumber \\
  U_{32} & = & cd - {\ts\half} h c^2, \nonumber \\
  U_{33} & = & d^2 + h\alpha cd - {\ts \frac{3}{4}} h^2 c^2.
        \label{eq:T-one-half-bis}
\end{eqnarray}
The elements of $T_{h,\alpha}^{1,\half}$ in~(\ref{eq:T-one-half})
and~(\ref{eq:T-one-half-bis}) may  be rewritten in terms of the alternative
choice
of deformation parameters~$h_{\pm}$.\par
%
%
The $T_h^{1,\half}$ matrix, corresponding to a single Jordanian deformation
parameter~$h$, may be obtained by choosing $\alpha=0$ (i.e. $h_+ = h_- = h$)
in~(\ref{eq:T-one-half}) and~(\ref{eq:T-one-half-bis}):
\begin{equation}
  T_h^{j=1} = D^{-1/2} \left(
      \begin{array}{ccc}
      a^2 + {\ts \frac{1}{4}} h^2 c^2 & \{a,b\} + {\ts \frac{1}{4}} h^2 \{c,d\}
           & 2b^2 + {\ts\half} h^2 \left(2D - 3a^2 + d^2\right) \\
                                                     &
           & - {\ts \frac{3}{8}} h^4 c^2 \\[0.2cm]
      {\ts\half} \{a,c\}                    & {\ts\half} \left(\{a,d\} +
\{b,c\}\right)
           & \{b,d\} - {\ts\frac{3}{4}} h^2 \{a,c\} \\[0.2cm]
      {\ts\half} c^2                          & {\ts\half} \{c,d\}
           & d^2 - {\ts \frac{3}{4}} h^2 c^2
      \end{array}\right).  \label{eq:T-one-h}
\end{equation}
Using the commutation relations, obtained from~(\ref{eq:Funh-alg}),
corresponding to the single-parametric case ($\alpha=0$), we have expressed the
elements of the $T_h^{j=1}$ matrix in~(\ref{eq:T-one-h}) in terms of the
determinant, and the anticommutators of the generating elements.\par
%
%
The contraction technique, comprising of the processes~(\ref{eq:tT})
and~(\ref{eq:T}), may be continued for obtaining the $(h,\alpha)$-Jordanian
monodromy matrices $T_{h,\alpha}^{j,z}$ for arbitrary coloured $(j,z)$
representations. Singularities systematically cancel leading to a well-defined
method of evaluation of the higher-dimensional $T_{h,\alpha}^{j,z}$
matrices.\par
%
%
\section{Conclusion}
\label{sec:conclusion}
Extending earlier results, we have shown that the $h$-deformed Jordanian
$R_h^{j_1;j_2}$ matrices may be obtained from the corresponding standard
$q$-deformed $R_q^{j_1;j_2}$ matrices through the present contraction
process. In
this demonstration, a key role is played by a nonlinear map of the U$_h$(sl(2))
algebra on the classical U(sl(2)) algebra. The Drinfeld twist operator
relating the
$h$-Jordanian and the classical coproduct structures may be determined as a
series expansion in~$h$ up to an arbitrary order. The twist operator
interrelating
the corresponding antipode maps has been determined in a closed form. By
using the
Reshetikhin formalism, the universal $\cR_{h,\alpha}$ matrix of the
two-parametric $(h,\alpha)$ Jordanian deformation has been obtained.\par
%
%
Moreover, we have shown that the two-parametric $(q,\lambda)$ standard
$R_{q,\lambda}^{j_1,z_1;j_2,z_2}$ matrices of arbitrary coloured $(j_1, z_1
\otimes j_2, z_2)$ representations reduce --- via the contraction procedure
--- to
the corresponding $(h,\alpha)$-Jordanian coloured matrices. \par
%
%
{}Finally, we have proved that the contraction process may serve as a
useful tool
to explicitly extract the so far unknown $(h,\alpha)$-Jordanian
$T_{h,\alpha}^{j,z}$ matrices for arbitrary $(j,z)$ coloured
representations from
the corresponding $(q,\lambda)$-deformed standard $T_{q,\lambda}^{j,z}$
matrices.\par
%
%
Our method may be useful in formulating Jordanian $h$-special functions in
analogy with the standard $q$-special functions. Of special relevance is the
well-known result~\cite{vaksman} that the elements of the standard $T_q^j$
matrices may be explicitly expressed in terms of the little $q$-Jacobi
polynomials.
It would be of interest to explore if a Jordanian analogue of this result
holds. We
hope to address this question elsewhere.
\par
%
%
\section*{Acknowledgments}
We thank A.~Chakrabarti, R.~Jagannathan, and J.~Segar for useful
discussions. One
of us~(RC) wishes to thank C.~Quesne for a kind invitation. He is also
grateful to
the members of the PNTPM~group for their kind hospitality.\par
%
%

\end{document}